\documentclass[a4paper,12pt]{amsart}

\usepackage[english]{babel}
\usepackage{amsthm}
 \usepackage{amsmath,amsfonts,amssymb}

\newcommand{\Set}[1]{\left\{\, #1 \,\right\}}

\DeclareMathOperator{\Sym}{Sym}

\DeclareMathOperator{\GL}{GL}
\DeclareMathOperator{\ML}{\mathfrak{M}}
\DeclareMathOperator{\Aff}{Aff}

\renewcommand{\phi}[0]{\varphi}
\renewcommand{\theta}[0]{\vartheta}
\renewcommand{\epsilon}[0]{\varepsilon}

\newcommand{\N}{\text{$\mathbf{N}$}}

\newcommand{\Q}{\text{$\mathbf{Q}$}}

\newcommand{\tran}[1]{\nu(#1)}
\newcommand{\TRAN}[1]{\tau(#1)}

\newtheorem{dummy}{Dummy}
\numberwithin{equation}{section}

\newtheorem{theorem}[dummy]{Theorem}
\newtheorem{lemma}[dummy]{Lemma}

\newtheorem{cor}[dummy]{Corollary}
\newtheorem{corollary}[dummy]{Corollary}

\theoremstyle{definition}

\theoremstyle{remark}

\newtheorem{fact}{Fact}
\newtheorem{example}{Example}

\newcommand{\oldplus}{+}
\newcommand{\oldminus}{-}
\newcommand{\newplus}{\circ}
\newcommand{\Newplus}{\diamond}

\begin{document}

\bibliographystyle{amsalpha}

\date{11 April 2006 --- Version 2.04%
}

\dedicatory{In memory of Edit Szab{\'o}}

\title[Abelian regular subgroups of the affine group]{Abelian regular subgroups of the affine group
and radical rings}

\author{A.~Caranti}

\address[A.~Caranti]{Dipartimento di Matematica\\
  Universit\`a degli Studi di Trento\\
  via Sommarive 14\\
  I-38050 Povo (Trento)\\
  Italy} 

\email{caranti@science.unitn.it} 

\urladdr{http://science.unitn.it/\~{ }caranti/}

\author{Francesca Dalla Volta}

\address[F.~Dalla Volta]{Dipartimento di Matematica e Applicazioni\\
  Edificio U5\\
  Universit\`a degli Studi di Milano--Bicocca\\
  Via R.~Cozzi, 53\\
  I-20126 Milano\\
  Italy}

\email{francesca.dallavolta@unimib.it}

\urladdr{http://www.matapp.unimib.it/\~{ }dallavolta/}

\author{Massimiliano Sala}

\address[M.~Sala]{Boole Centre for Research in Informatics\\
University College Cork\\
Cork\\
Ireland}

\email{msala@bcri.ucc.ie}

\begin{abstract}
   We establish  a correspondence between abelian  regular subgroup of
   the affine  group, and commutative,  associative algebra structures
   on the underlying vector space that are (Jacobson) radical rings.

   As  an  application, we  show  that  if  the underlying  field  has
   positive  characteristic,  then  an  abelian regular  subgroup  has
   finite exponent if the vector space is finite-dimensional, while it
   can be torsion free if the dimension is infinite.

   We  also give an  example of  an abelian,  regular subgroup  of the
   affine  group  over  an  infinite-dimensional vector  space,  which
   intersects trivially the group of translations.
\end{abstract}

\keywords{affine group, abelian regular subgroups, (Jacobson) radical rings}

\subjclass{20B10 16N20}

\thanks{First author partially supported by MIUR-Italy via PRIN 2003018059
 ``Graded Lie algebras and pro-$p$-groups: representations,
 periodicity and derivations''.  Second  author partially  supported by
 MIUR-Italy via PRIN ``Group theory and applications''. First two
 authors are members of INdAM-GNSAGA}

\maketitle

\thispagestyle{empty}

\section{Introduction}

Cai  Heng   Li  has   described  in~\cite{Li}  the   finite  primitive
permutation groups which contain an abelian regular subgroup. Among these
we have the groups of affine  type, where the translations form an abelian
regular subgroup  which is normal. Li  notes that there  might well be
other (non normal) abelian regular subgroups in such groups, and as an
example  describes the  abelian  regular subgroups  of  the group  of
affine  type  that  can be  obtained  as  the  split extension  of  an
elementary abelian  group of order  $2^{d}$ by the symmetric  group on
$d+1$ letters.

The goal of  this  note is to record a  simple description of the
abelian regular  subgroups of the full  affine group in  terms of
commutative, associative algebra structures that one can impose on the
underlying vector space, so that the resulting ring is radical.

In   Section~\ref{sec:regular}   we   establish   the   correspondence
(Theorem~\ref{theorem:general}). In Section~\ref{sec:examples} we give
some examples.

As an  application, we prove  in Corollary~\ref{cor:appl} that  if the
underlying field has positive  characteristic, then an abelian regular
subgroup has finite exponent if the vector space is
finite-dimensional. In Example~\ref{ex:general} we show that if the
dimension is allowed to be infinite, then an abelian regular subgroup
can be torsion free.

An abelian  regular subgroup of an  affine group over  a finite vector
space must intersect the  group of translations nontrivially.  P{\'a}l
Heged{\H{u}}s  has given  an example~\cite{Hegedus}  of  a nonabelian,
regular subgroup of  an affine group over a  finite vector space which
has  trivial   intersection  with   the  group  of   translations.  In
Corollary~\ref{theorem:example}  we  show   that  the  same  setting  of
Example~\ref{ex:general}  provides an example  of an  abelian, regular
subgroup of the affine group over an infinite-dimensional vector space
which has trivial intersection with the group of translations.

\section{Abelian regular subgroups}
\label{sec:regular}

Let $F$ be an arbitrary field, and
$(V, \oldplus)$ be a
vector  space  of arbitrary  dimension $d$  over $F$.

Let $\GL(V)$ be the group of invertible, $F$-linear
maps on $V$, and $N$ be the group of translations, that is,
\begin{equation*}
  N = \Set { \tran{x} : x \in V },
\end{equation*}
where $\tran{x} : z \mapsto z \oldplus x$.
Let  $\Aff(V)$  be  the    affine
group on $V$.

Clearly we have
\begin{fact}\label{eq:afffact}
$N$ is a normal subgroup of $\Aff(V)$. Every element of $\Aff(V)$ can
 be written \emph{uniquely} as a product of an element of $\GL(V)$ and an
 element of $N$,
 so that we have the semidirect product decomposition
 \begin{equation*}
   \Aff(V) = \GL(V) N.
 \end{equation*}
\end{fact}
We also write $\ML(V)$ for the $F$-algebra of $F$-linear
maps on $V$.

Recall that a  group $G$ of permutations on a set  $\Omega$ is said to
be \emph{regular}  if, given  any $\alpha \in  \Omega$, then  for each
$\beta \in \Omega$ there exists a unique $g \in G$ such that $\alpha g
= \beta$.

Clearly $N$ is an abelian regular subgroup of $\Aff(V)$. The next
result describes them all.

Recall  that   a  \emph{(Jacobson)  radical  ring}~\cite[Definition~2,
p.~4]{Jac} is a ring $(A,  \oldplus, \cdot)$ in which every element is
invertible with respect to the circle  operation $x \newplus y = x + y
+ x \cdot y$, so that $(A, \newplus)$ is a group. Equivalently, a ring
is radical if it coincides with its Jacobson radical.

\begin{theorem}
  \label{theorem:general}
  Let $F$ be an arbitrary field, and $(V, \oldplus)$ a 
  vector space of arbitrary dimension over $F$.

  There is a one-to-one correspondence between
  \begin{enumerate}
  \item abelian regular subgroups $T$ of $\Aff(V)$, and
  \item commutative, associative $F$-algebra structures $(V,
    \oldplus, \cdot)$ that one can impose on the vector space structure
    $(V, \oldplus)$, such that the resulting ring is radical.
  \end{enumerate}

  In   this  correspondence,   isomorphism  classes   of  $F$-algebras
  correspond to conjugacy  classes under the action of  $\GL(V)$ of
  abelian regular subgroups of  $\Aff(V)$.
\end{theorem}

\begin{proof}

Let $T$  be an  abelian regular subgroup  of $\Aff(V)$.  Since  $T$ is
regular, for each $x \in V$ there is a unique $\TRAN{x} \in T$ such that $0
\TRAN{x} = x$. (Our affine maps act on the right.) Thus
\begin{equation}\label{eq:T}
  T = \Set { \TRAN{x} : x \in V }.
\end{equation}
Because of Fact~\ref{eq:afffact}, we can write $\TRAN{x}$ uniquely as
\begin{equation}\label{eq:fact}
\TRAN{x} = \gamma(x) \tran{x},
\end{equation}
where $\gamma(x) \in \GL(V)$. We also introduce
\begin{equation*}
  \delta(x) = \gamma(x) \oldminus 1 \in \ML(V).
\end{equation*}

We have
\begin{equation}\label{eq:product}
\begin{aligned}
  \TRAN{x} \TRAN{y}
  &=
  \gamma(x) \tran{x} \gamma(y) \tran{y}
  \\&=
  \gamma(x) \gamma(y) \tran{x}^{\gamma(y)} \tran{y}
  \\&=
  \gamma(x) \gamma(y) \tran{x \gamma(y)} \tran{y}
  \\&=
  \gamma(x) \gamma(y) \tran{x \oldplus y \oldplus x \delta(y)}.
\end{aligned}
\end{equation}
As $T$ is abelian, we have $\TRAN{x} \TRAN{y} = \TRAN{y} \TRAN{x}$
for all $x, y \in V$. 
Therefore~\eqref{eq:product} yields
\begin{equation*}
  \gamma(x) \gamma(y) \tran{x \oldplus y \oldplus x \delta(y)}
  =
  \gamma(y) \gamma(x) \tran{y \oldplus x \oldplus y \delta(x)}
\end{equation*}
for all $x, y \in V$. 
From Fact~\ref{eq:afffact} we get
\begin{equation*}
  \tran{x \oldplus y \oldplus x \delta(y)}
  =
  \tran{y \oldplus x \oldplus y \delta(x)}
\end{equation*}
for all $x, y \in V$. 
We get
\begin{fact}
  \label{eq:commutes} $x \delta(y) = y \delta(x)$ for all $x, y \in
  V$.
\end{fact}
As the left-hand side of $x \delta(y) = y \delta(x)$ is linear in
$x$, so is the right-hand side. We obtain
\begin{fact}
  \label{eq:distributes}
  $\delta : V \to \ML(V)$ is $F$-linear.
\end{fact}

Since $T$ is a group, we have $\TRAN{x} \TRAN{y} = \TRAN{z}$ for some $z
\in V$.  Because of~\eqref{eq:product} and  Fact~\ref{eq:afffact}, we
have $z = x \oldplus y \oldplus x \delta(y)$, so that
\begin{fact}
  \label{fact}
  $\TRAN{x} \TRAN{y} = \TRAN{ x \oldplus y \oldplus x \delta(y)}$ for
  all $x, y \in V$.
\end{fact}
Thus,  again  by~\eqref{eq:product}  and  Fact~\ref{eq:afffact},
$\gamma(x) \gamma(y)  = \gamma(x \oldplus y \oldplus  x \delta(y))$. We
obtain, using Fact~\ref{eq:distributes},
\begin{align*}
  \gamma(x) \gamma(y)
  &=
  1 \oldplus \delta(x) \oldplus \delta(y) \oldplus \delta(x) \delta(y)
  \\&=
  \gamma(x \oldplus y \oldplus x \delta(y))
  \\&=
  1 \oldplus \delta(x \oldplus y \oldplus x \delta(y))
  \\&=
  1 \oldplus \delta(x) \oldplus \delta(y) \oldplus \delta(x \delta(y)),
\end{align*}
that is, 
\begin{fact}\label{eq:preass}
  $\delta(x \delta(y)) = \delta(x) \delta(y)$ for all $x, y \in V$.
\end{fact}

Now define on $V$ a product operation by
\begin{equation}
  \label{eq:cdot}
  x \cdot y = x \delta(y).
\end{equation}
This  product  is   commutative,  by  Fact~\ref{eq:commutes}.   It  is
$F$-linear  in  both  variables  (in particular  it  distributes  over
$\oldplus$),  for  instance  because  $\delta(y)  \in  \ML(V)$  is  an
$F$-linear map, and because of Fact~\ref{eq:distributes}.  The product
is  also  associative  as  for  all  $x,  y, z  \in  V$  one  has,  by
definition~\eqref{eq:cdot} and Fact~\ref{eq:preass},
\begin{equation*}
 (x y) z = (x \delta(y)) z = x \delta(y) \delta(z)
 =
 x \delta(y \delta(z)) = x \delta(y z) = x (y z).
\end{equation*}
Therefore $(V, +, \cdot)$ is an $F$-algebra.

We can now consider the circle operation ``$\newplus$'' on $V$ given by
\begin{equation}
  \label{eq:circle}
  x \newplus y = x \oldplus y \oldplus x y,
\end{equation}
which makes $(V, \newplus)$ into a monoid. The map
\begin{align*}
\tau: (V, \newplus) & \to T\\
      x             & \mapsto \TRAN{x}
\end{align*}
is an isomorphism of monoids,
because    Fact~\ref{fact}   can    be    rewritten,   according    to
Fact~\ref{eq:preass},    \eqref{eq:cdot}   and~\eqref{eq:circle},   as
$\tau(x) \tau(y) = \tau(x \newplus y)$.
Since $T$ is a group, so is $(V, \newplus)$. We have obtained
\begin{fact}\label{part:amap}
  $(V, \newplus)$ is an abelian group, and the map
  \begin{align*}
    \tau: (V, \newplus) & \to T\\
    x             & \mapsto \TRAN{x}
  \end{align*}
  is a group isomorphism.
\end{fact}
We have thus proved that the ring $(V, +, \cdot)$ is radical.

We note also the following
\begin{fact}
  \label{eq:TRAN}
  $z \TRAN{x} = z \newplus x$.
\end{fact}
This follows from
\begin{equation*}
  z \TRAN{x} = z (1 \oldplus \delta(x)) \tran{x}
  = z \oldplus z \delta(x) \oldplus x = z \newplus x.
\end{equation*}
Fact~\ref{eq:TRAN} shows  that $\TRAN{x}$
is  also a  translation,  but with  respect  to ``$\newplus$'',  while
$\tran{x}$ is a translation  with respect to ``$\oldplus$''. (But note
that $(V, \newplus)$ need not be the additive group of a vector space,
see e.g.\ Example~\ref{ex:S4} in the next section.)

Conversely, suppose $(V,  +, \cdot)$ is a radical ring.  For $x \in V$
define a map  $\tau(x)$ on $V$ by $\tau(x) : y  \mapsto y \newplus x$,
with  ``$\newplus$''  as  in~\eqref{eq:circle}.  Reversing  the  above
arguments, one  sees that  $T = \Set  { \TRAN{x}  : x \in  V }$  is an
abelian regular subgroup of $\Aff(V)$ then.

  We now  pass to  the statement about  the isomorphism  and conjugacy
  classes.

  Suppose first that $(V, \oldplus, \cdot)$ and $(V, \oldplus, *)$ are two
  commutative, associative $F$-algebra structures on the vector space
  structure $(V, +)$, such that they are radical rings. Suppose
  there is an $F$-algebra isomorphism 
  \begin{equation*}
    \phi : (V, \oldplus, \cdot) \to (V, \oldplus, *).
  \end{equation*}
  In particular, $\phi \in \GL(V)$. For $x, y \in V$ we have
  two circle operations
  \begin{equation*}
    \begin{cases}
      x \newplus y = x \oldplus y \oldplus x \cdot y,\\
      x \Newplus y = x \oldplus y \oldplus x * y,\\
    \end{cases}
  \end{equation*}
  Since $\phi$  is  an algebra  isomorphism it  follows
  \begin{equation}
    \label{eq:mor}
    \begin{aligned}
    (x \newplus y) \phi 
    &=
    (x \oldplus y \oldplus x \cdot y) \phi
    \\&=
    (x \phi) \oldplus (y \phi) \oldplus (x \phi) * (y \phi)
    \\&= 
    x \phi  \Newplus y \phi.
  \end{aligned}
    \end{equation}
  
  Let $T_{1} = \Set
  { \tau_{1}(x) : x  \in V }$ and $T_{2} = \Set  { \tau_{2}(x) : x \in
  V}$ be the corresponding subgroups, so that
  \begin{equation*}
    \begin{cases}
      z \tau_{1}(x) = z \newplus x,\\
      z \tau_{2}(x) = z \Newplus x.\\
    \end{cases}
  \end{equation*}

  Now for $x, y \in  V$ we have, according to Fact~\ref{eq:TRAN} and
  to~\eqref{eq:mor},
  \begin{equation*}
     y \phi^{-1} \tau_{1}(x) \phi
     =
     (y \phi^{-1} \newplus x) \phi
     =
     y \Newplus x \phi
     =
     y \tau_{2}( x \phi ).
  \end{equation*}
  Thus for all $x \in V$ we have $\tau_{2}(x \phi) =
  \phi^{-1}  \tau_{1}(x) \phi $, so that $T_{2} =  \phi^{-1} T_{1} \phi$.

  Conversely, if $T_{2}  = \phi^{-1} T_{1} \phi$ for some $\phi \in
  \GL(V)$, let $\psi : V \to V$ be the bijection such that
  \begin{equation*}
    \phi^{-1} \tau_{1}(x) \phi = \tau_{2} (x \psi)
  \end{equation*}
  for all $x \in V$. We have
  \begin{equation*}
    0 \phi^{-1} \tau_{1}(x) \phi 
    = 
    0 \tau_{1}(x) \phi 
    =
    x \phi 
    =
    0 \tau_{2} (x \psi)
    =
    x \psi,
  \end{equation*}
  so that $\phi = \psi$.
  
  Now for $x, y \in V$ we have, according to Fact~\ref{part:amap}, 
  \begin{align*}
    \tau_{2} ( x \phi \Newplus y \phi)
    &=
    \tau_{2}(x \phi) \tau_{2}(y \phi)
    \\&=
    \phi^{-1} \tau_{1}(x) \tau_{1}(y) \phi
    \\&=
    \phi^{-1} \tau_{1}(x \newplus y) \phi
    \\&=
    \tau_{2} ( (x \newplus y) \phi),
  \end{align*}
  so that 
  \begin{equation*}
    (x \newplus y) \phi = x \phi \Newplus y \phi,
  \end{equation*}
  and reversing the argument of~\eqref{eq:mor} one gets that
  \begin{equation*}
    \phi : (V, \oldplus, \cdot) \to (V, \oldplus, *)
  \end{equation*}
  is an isomorphism of $F$-algebras.
\end{proof}

As an application, we get
\begin{cor}
  \label{cor:appl}
  Let $F$ be a field of positive characteristic $p$, and let $(V,
  \oldplus)$ be a finite dimensional vector space over $F$.

  Then every abelian regular subgroup of $\Aff(V)$ has finite
  exponent, which is a power of $p$.
\end{cor}

The result does not hold when $V$ is allowed to be infinite dimensional. In
Example~\ref{ex:general} in the next section we give an example of an
abelian regular subgroup which is torsion free.

\begin{proof}
  We first prove that any $F$-algebra $(V, \oldplus, \cdot)$, which is
  radical as a ring, is  nilpotent. Clearly $V$ is non-unital, because
  $-1$  would have  no inverse  with respect  to $\newplus$,  as $(-1)
  \newplus  a   =  -1$   for  all  $a   \in  V$.   Now   a  non-unital
  finite-dimensional $F$-algebra $A$ need  not be an Artinian ring, as
  an ideal need not be an $F$-subspace.  (Think of the one-dimensional
  $\Q$-algebra $\Q  a$, where $a^{2}  = 0$.  Its ideals  correspond to
  the additive  subgroups of $\Q$.)   However, $A$ can be  embedded in
  the standard way in a  unital $F$-algebra $B$ of dimension one more.
  If $A$ is radical, then $A$  is the Jacobson radical of the artinian
  ring $B$, and so $A$ is nilpotent.

  For $a \in \N$ and $x \in
  V$, we write 
  \begin{equation*}
    a_{\newplus} x 
    = 
    \underbrace{x \newplus \dots \newplus x}_{\text{$a$ times}}.
  \end{equation*}
  One proves easily by induction
  \begin{equation*}
    a_{\newplus} x = \sum_{i = 1}^{a} \binom{a}{i} x^{i}.
  \end{equation*}
  In particular
  \begin{equation}
    \label{eq:exponents}
    p^{j}_{\newplus} x = x^{p^{j}}.
  \end{equation}
  
  Suppose $V^{n}  = 0$ for  some $n$. If  $p^{j} \ge n$ for  some $j$,
  then   $V^{p^{j}}  =   0$,  so   that   by~\eqref{eq:exponents},  and
  Fact~\ref{part:amap}, the corresponding
  $T$ has exponent dividing $p^{j}$.
\end{proof}

\section{Examples and comments}
\label{sec:examples}

We  begin with  some examples  in  which $V$  is finite,  so that  all
algebra structures $(V, \oldplus, \cdot)$, which are radical as rings,
are nilpotent. In this section we use the notation of the proof of
Theorem~\ref{theorem:general}. 

The first example is here as folklore.
\begin{example}
  \label{ex:S4}
  When $F$ is the field with $2$ elements, and $V$ has dimension $2$
  over $F$, the affine group is isomorphic to
  the symmetric group $\Sym(4)$ on four letters. We write $V = \Set{0,
  a, b, a \oldplus b}$ for the underlying vector space. The other
  abelian regular subgroups of $\Sym(4)$ are the three cyclic
  subgroups, which correspond to the ring structures on $(V, \oldplus)$
  defined by 
  \begin{itemize}
  \item $a^{2} = b$, $b^{2} = 0$, $a b = 0$;
  \item $a^{2} = 0$, $b^{2} = a$, $a b = 0$;
  \item $a^{2} = b^{2} = a b = a \oldplus b$.
  \end{itemize}
  For instance in the first case we obtain the cyclic group $(T, \newplus)$
  where 
  \begin{itemize}
  \item $2_{\newplus}  a =  a \oldplus a  \oldplus b =  b$, 
  \item $3_{\newplus}  a = b \oldplus a  \oldplus a b =  a \oldplus b$, 
  \item $4_{\newplus}  a = a \oldplus  b \oldplus a \oldplus  (a
  \oldplus b)  a = b \oldplus b = 0$.
  \end{itemize}
  The three cyclic  subgroups are conjugate here, and  the three rings
  are isomorphic.
\end{example}
A generalization  of this for an  arbitrary prime $p$ is  given by the
following
\begin{example}
  Let $p$ be an arbitrary prime, and take $V$ to have dimension $p$
  over the field $F$ with $p$ elements. Then
  one can define a suitable ring structure on $V$ by declaring a base of $V$ in
  the form
  \begin{equation*}
    a, a^{2}, \dots, a^{p},
  \end{equation*}
  and letting $a^{p+1} = 0$. The corresponding group $T$ is abelian of
  type $(p^{2}, \underbrace{p, \dots, p}_{p-2})$, where the cyclic
  component of order $p^{2}$ is generated by $a$ (one has
  $p_{\newplus} a = a^{p}$), and those of order $p$ are generated by $a^{2},
  \dots, a^{p-1}$.
\end{example}

The following non trivial ring structure on $V$ satisfies $x y z = 0$
for all $x, y, z \in V$.
\begin{example}
  \label{ex:normalizes}
  Let  $F$ be  the field  with $2$  elements, and   $(V, \oldplus,
  \cdot)$ be  the exterior  algebra over a  vector space  of dimension
  $k$, spanned by $e_{1}, \dots, e_{k}$, truncated at length $2$. That
  is, $V$ has basis
  \begin{equation*}
    e_{1}, \dots, e_{k}, e_{1} \wedge e_{2}, \dots , e_{k-1} \wedge e_{k},
  \end{equation*}
  and satisfies $x^{2} = 0$.
\end{example}
This is  relevant to the  question whether $N$ normalizes  all abelian
regular subgroups $T$. Note  first the following interpretation of our
product in terms of the action of $T$ on $N$.
\begin{lemma}\label{lemma:comm}
  $[\tran{x}, \TRAN{y}] = \tran{x y}$.
\end{lemma}
\begin{proof}
  \begin{align*}
    [\tran{x}, \TRAN{y}] 
    &= \tran{x}^{-1} \tran{x}^{\gamma(y)} 
    \\&= \tran{- x}  \tran{x (1 \oldplus \delta(y))} 
    \\& = \tran{- x \oldplus x \oldplus x y}  = \tran{x y}.
  \end{align*}
\end{proof}
By assumption, $T$ normalizes $N$. Now $N$ normalizes $T$ if
and only if $\tran{x y} \in T$ for all $x, y$, that is, $\delta(x y) =
0$, that is $x y z = 0$ for all $x, y, z \in V$. So in
Example~\ref{ex:normalizes}, $N$ normalizes $T$. However, in the
following example $N$ does not normalize $T$, as there is a nonzero
threefold product.
\begin{example}
  \label{ex:exterior}
  Let $F$ be  the field with $2$ elements,  and $(V, \oldplus, \cdot)$
  be  the exterior  algebra over  a vector  space of  dimension three,
  spanned by $e_{1}, e_{2}, e_{3}$. That is, $V$ has basis
\begin{equation*}
  e_{1}, e_{2}, e_{3}, e_{1} \wedge e_{2}, e_{1} \wedge e_{3} , e_{2}
  \wedge e_{3}, e_{1} \wedge e_{2} \wedge e_{3}.
\end{equation*}
Clearly $x^{2} = 0$ for all $x$, but $e_{1} \wedge e_{2} \wedge e_{3}
\ne 0$. 
\end{example}

In Examples~\ref{ex:normalizes} and~\ref{ex:exterior},  the ring is an
exterior algebra over  the field $F$ with two  elements, or a quotient
thereof.  In  characteristic  $2$,  algebras  that  are  quotients  of
exterior algebras  correspond to elementary  abelian regular subgroups
of the affine group.

Before considering an example when $V$ is infinite, let us define, for
a prescribed $F$-algebra structure $(V, +, \cdot)$,
 \begin{equation*}
   U 
   = 
   \ker(\delta) 
   = 
   \Set{x \in V : \text{$x \cdot y = 0$ for all $y \in   V$}}.
 \end{equation*}
Clearly we can choose the  algebra structure on the finite dimensional
vector space $(V,  \oldplus)$ so that $U$ has  arbitrary dimension, for
instance as in the next example.
\begin{example}
  \label{ex:e+1}
  Let $0  \le e < d$, and we choose  $(V, \oldplus, \cdot)$ to
  be the quotient  of the ideal of the  polynomial ring $F[x_{0}, x_{1},
    \dots, x_{e}]$  generated by $x_{0}, x_{1}, \dots,  x_{e}$ modulo the
  relations
  \begin{equation*}
    \begin{cases}
      x_{0}^{d - e + 1},\\
      x_{0} x_{i}, & \text{for $i > 0$},\\
      x_{j} x_{i}, & \text{for $i, j > 0$}.\\
    \end{cases}
  \end{equation*}
  Then $V$ has dimension $d$, while $U$ has basis $ x_{0}^{d - e},
  x_{1},  \dots, x_{e}$, and thus $\dim(U) = e + 1$.
\end{example}

Now $U$ corresponds to the intersection $N \cap T$, as
\begin{equation}
  \label{eq:NcapT-U}
  \begin{aligned}
    N \cap T 
    &= 
    \Set{ \tran{x} : \TRAN{x} = \tran{x}}
    \\&= 
    \Set{ \tran{x} : \delta(x) = 0} 
    \\&= 
    \Set{ \tran{x} : x \in U}.
  \end{aligned}
\end{equation}
Using Lemma~\ref{lemma:comm}, one recovers the well-known fact
\begin{lemma}
  Let $N$ be the group  of translations in the affine group $\Aff(V)$,
  and  let $T$  be  an abelian  regular  subgroup. Then 
  \begin{equation}
    \label{eq:NcapT}
    N  \cap T  =  C_{N}(T) = C_{T}(N).
  \end{equation}
\end{lemma}

When $V$, and thus
$\Aff(V)$, is  finite, then $U \ne 0$, as $(V,  +, \cdot)$ is nilpotent,
so  that  the subgroup  of~\eqref{eq:NcapT}  is
nontrivial. (Alternatively, $C_{N}(T)$ is nontrivial, as $T$ is a
finite $p$-group acting on the finite $p$-group $N$; here $p$ is the
characteristic of the underlying field.) In other
words, an abelian regular subgroup of the affine group over  a finite
vector space intersects the group of translations nontrivially.

It  also follows from  Example~\ref{ex:e+1} that  when $V$  is finite,
then $N \cap T$ has arbitrary order, different from $1$.

P{\'a}l Heged{\H{u}}s has given an example~\cite{Hegedus} of a
nonabelian, regular subgroup of an affine group over a finite vector
space which has trivial intersection with the group of
translations. 

Now we consider the following
\begin{example}
  \label{ex:general}
  Let $(V,  \oldplus, \cdot)$ be the  maximal ideal $t  F[[t]]$ of the
  $F$-algebra $F[[t]]$ of formal  power series over an arbitrary field
  $F$. This is a radical ring.  Since $F[[t]]$ is a domain, we have $U
  = 0$ here.

  It follows from~\eqref{eq:NcapT-U} that in this example the abelian
  regular subgroup $T$ intersects trivially the group $N$ of translations.

  Also, $T$ is torsion-free. If $F$ is a field of positive
  characteristic $p$, then the group $N$ of translations has exponent
  $p$. Thus $\Aff(V)$ has two rather different abelian
  regular subgroups here.
\end{example}

Summing up, we have
\begin{corollary}
  \label{theorem:example}
  \

  \begin{enumerate}
    \item In the  affine group over a finite  vector space, an abelian
    regular   subgroup   intersects    the   group   of   translations
    nontrivially.
  \item There  is an  example~\cite{Hegedus} of a  nonabelian, regular
    subgroup of an  affine group over a finite  vector space which has
    trivial intersection with the group of translations.
  \item There is an example (Example~\protect{\ref{ex:general}} above)
    of  an abelian,  regular  subgroup  of the  affine  group over  an
    infinite  vector space  which  has trivial  intersection with  the
    group of translations.
  \end{enumerate}
\end{corollary}

\providecommand{\bysame}{\leavevmode\hbox to3em{\hrulefill}\thinspace}
\providecommand{\MR}{\relax\ifhmode\unskip\space\fi MR }
\providecommand{\MRhref}[2]{%
  \href{http://www.ams.org/mathscinet-getitem?mr=#1}{#2}
}
\providecommand{\href}[2]{#2}

\end{document}